\input amstex
\let\myfrac=\frac
\input eplain.tex
\let\frac=\myfrac
\input epsf




\loadeufm \loadmsam \loadmsbm
\message{symbol names}\UseAMSsymbols\message{,}

\font\myfontdefault=cmr10

\font\mytdmchapfont=cmb10 at 14pt
\font\mytdmheadfont=cmb10 at 10pt
\font\mytdmsubheadfont=cmr10

\magnification 1200
\newif\ifinappendices
\newif\ifundefinedreferences
\newif\ifchangedreferences
\newif\ifloadreferences
\newif\ifmakebiblio
\newif\ifmaketdm

\undefinedreferencesfalse
\changedreferencesfalse


\loadreferencestrue
\makebibliofalse
\maketdmfalse

\def\headpenalty{-400}     
\def\proclaimpenalty{-200} 

%
%

\def\alphanum#1{\ifcase #1 _\or A\or B\or C\or D\or E\or F\or G\or H\or I\or J\or K\or L\or M\or N\or O\or P\or Q\or R\or S\or T\or U\or V\or W\or X\or Y\or Z\fi}
\def\gobbleeight#1#2#3#4#5#6#7#8{}

\newwrite\references
\newwrite\tdm
\newwrite\biblio

\newcount\chapno
\newcount\headno
\newcount\subheadno
\newcount\procno
\newcount\figno
\newcount\citationno

\def\setcatcodes{%
\catcode`\!=0 \catcode`\\=11}%

\ifloadreferences
    {\catcode`\@=11 \catcode`\_=11%
    \global\def\_@citation@BallGromSch{1}
\global\def\_@citation@Grom{2}
\global\def\_@citation@LabC{3}
\global\def\_@citation@LabB{4}
\global\def\_@citation@LabA{5}
\global\def\_@citation@RosSpruck{6}
\global\def\_@citation@SmiB{7}
\global\def\_@citation@SmiC{8}
\global\def\_@citation@SmiE{9}
\global\def\_@proc@Area_Is_Finite{1.1}
\global\def\_@proc@Volume_Is_Finite{1.2}
\global\def\_@head@Polar_Coordinates_About_A_Geodesic{2}
\global\def\_@proc@Metric_In_Polar_Coordinates{2.2}
\global\def\_@proc@Connection_One_Form_In_Polar_Coordinates{2.3}
\global\def\_@head@Surfaces_Of_Revolution{3}
\global\def\_@proc@Gaussian_Curvature_Of_Surface_Of_Revolution{3.1}
\global\def\_@head@Surfaces_Of_Constant_Curvature{4}
\global\def\_@proc@The_Function_Decays_Exponentially{4.1}
\global\def\_@proc@Control_Of_Area_And_Volume_Of_Rotationally_Symmetric_Cusp{4.3}
\global\def\_@head@Asymptotically_Tubular_Immersed_Surfaces{5}
\global\def\_@proc@Surface_Is_Transverse{5.1}
\global\def\_@proc@Intersection_Is_Smooth_Curve{5.2}
\global\def\_@head@Convex_Curves_In_Real_And_Hyperbolic_Space{6}
\global\def\_@proc@Winding_Number_Is_Bounded_By_Index{6.1}
\global\def\_@proc@Stokes_Theorem{6.2}
\global\def\_@proc@Derivative_Of_Winding_Number{6.3}
\global\def\_@head@The_Volume_Contained_By_A_Cusp{7}
\global\def\_@proc@Image_Lies_Inside_Cusp{7.1}
\global\def\_@proc@V_Is_Volume{7.2}
\global\def\_@proc@The_Volume_Function_Is_Bounded{7.3}
\global\def\_@head@Finiteness_Of_Area_And_Volume{8}
    }%
\else
    \openout\references=references.tex
\fi

\newcount\newchapflag 
\newcount\showpagenumflag 

\global\chapno = -1 
\global\citationno=0
\global\headno = 0
\global\subheadno = 0
\global\procno = 0
\global\figno = 0

\def\resetcounters{%
\global\headno = 0%
\global\subheadno = 0%
\global\procno = 0%
\global\figno = 0%
}

\global\newchapflag=0 
\global\showpagenumflag=0 

\def\chinfo{\ifinappendices\alphanum\chapno\else\the\chapno\fi}%
\def\headinfo{\the\headno}
\def\subheadinfo{\the\headno.\the\subheadno}
\def\procinfo{\the\headno.\the\procno}
\def\figinfo{\the\headno.\the\figno}
\def\citationinfo{\the\citationno}%
\def\nextheadno{\global\advance\headno by 1 \global\subheadno = 0 \global\procno = 0}
\def\nextsubheadno{\global\advance\subheadno by 1}
\def\nextprocno{\global\advance\procno by 1 \procinfo}
\def\nextfigno{\global\advance\figno by 1 \figinfo}

{\global\let\noe=\noexpand%
%
%
\catcode`\@=11%
\catcode`\_=11%
\setcatcodes%
!global!def!_@@internal@@makeref#1{%
!global!expandafter!def!csname #1ref!endcsname##1{%
!csname _@#1@##1!endcsname%
!expandafter!ifx!csname _@#1@##1!endcsname!relax%
    !write16{#1 ##1 not defined - run saving references}%
    !undefinedreferencestrue%
!fi}}%
!global!def!_@@internal@@makelabel#1{%
!global!expandafter!def!csname #1label!endcsname##1{%
!edef!temptoken{!csname #1info!endcsname}%
!ifloadreferences%
    !expandafter!ifx!csname _@#1@##1!endcsname!relax%
        !write16{#1 ##1 not hitherto defined - rerun saving references}%
        !changedreferencestrue%
    !else%
        !expandafter!ifx!csname _@#1@##1!endcsname!temptoken%
        !else
            !write16{#1 ##1 reference has changed - rerun saving references}%
            !changedreferencestrue%
        !fi%
    !fi%
!else%
    !expandafter!edef!csname _@#1@##1!endcsname{!temptoken}%
    !edef!textoutput{!write!references{\global\def\_@#1@##1{!temptoken}}}%
    !textoutput%
!fi}}%
!global!def!makecounter#1{!_@@internal@@makelabel{#1}!_@@internal@@makeref{#1}}%
!unsetcatcodes%
}
\makecounter{ch}%
\makecounter{head}%
\makecounter{subhead}%
\makecounter{proc}%
\makecounter{fig}%
\makecounter{citation}%
\def\newref#1#2{%
\def\temptext{#2}%
\edef\bibliotextoutput{\expandafter\gobbleeight\meaning\temptext}%
\global\advance\citationno by 1\citationlabel{#1}%
\ifmakebiblio%
    \edef\fileoutput{\write\biblio{\noindent\hbox to 0pt{\hss$[\the\citationno]$}\hskip 0.2em\bibliotextoutput\medskip}}%
    \fileoutput%
\fi}%
\def\cite#1{%
$[\citationref{#1}]$%
\ifmakebiblio%
    \edef\fileoutput{\write\biblio{#1}}%
    \fileoutput%
\fi%
}%
%
%
%

\let\mypar=\par


\def\raggedleft{\leftskip=0pt plus 1fil \parfillskip=0pt}


\font\lettrinefont=cmr10 at 28pt
\def\lettrine #1[#2][#3]#4%
{\hangafter -#1 \hangindent #2
\noindent\hskip -#2 \vtop to 0pt{
\kern #3 \hbox to #2 {\lettrinefont #4\hss}\vss}}

\font\mylettrinefont=cmr10 at 28pt
\def\mylettrine #1[#2][#3][#4]#5%
{\hangafter -#1 \hangindent #2
\noindent\hskip -#2 \vtop to 0pt{
\kern #3 \hbox to #2 {\mylettrinefont #5\hss}\vss}}


\edef\Pagetitle={Blank}

\headline={\hfil\Pagetitle\hfil}

\footline={{\myfontdefault \hfil\folio\hfil}}

\def\nextoddpage
{
\newpage%
\ifodd\pageno%
\else%
    \global\showpagenumflag = 0%
    \null%
    \vfil%
    \eject%
    \global\showpagenumflag = 1%
\fi%
}


\def\newchap#1#2%
{%
%
%
\global\advance\chapno by 1%
\resetcounters%
%
%
\newpage%
\ifodd\pageno%
\else%
    \global\showpagenumflag = 0%
    \null%
    \vfil%
    \eject%
    \global\showpagenumflag = 1%
\fi%
\global\newchapflag = 1%
\global\showpagenumflag = 1%
%
%
{\font\chapfontA=cmsl10 at 30pt%
\font\chapfontB=cmsl10 at 25pt%
\null\vskip 5cm%
{\chapfontA\raggedleft\hfil%
{%
\ifnum\chapno=0
    \phantom{%
    \ifinappendices%
        Annexe \alphanum\chapno%
    \else%
        \the\chapno%
    \fi}%
\else%
    \ifinappendices%
        Annexe \alphanum\chapno%
    \else%
        \the\chapno%
    \fi%
\fi%
}%
\par}%
\vskip 2cm%
{\chapfontB\raggedleft%
\lineskiplimit=0pt%
\lineskip=0.8ex%
\hfil #1\par}%
\vskip 2cm%
}%
\edef\Pagetitle{#2}%
%
%
\ifmaketdm%
    \def\temp{#2}%
    \def\tempbis{\nobreak}%
    \edef\chaptitle{\expandafter\gobbleeight\meaning\temp}%
    \edef\mynobreak{\expandafter\gobbleeight\meaning\tempbis}%
    \edef\textoutput{\write\tdm{\bigskip{\noexpand\mytdmchapfont\noindent\chinfo\ - \chaptitle\hfill\noexpand\folio}\par\mynobreak}}%
\fi%
\textoutput%
}


\def\newhead#1%
{%
\ifhmode%
    \mypar%
\fi%
\ifnum\headno=0%
\else%
   \nobreak\vskip -\lastskip%
   \penalty\headpenalty\vskip .5cm%
\fi%
\nextheadno%
\ifmaketdm%
    \def\temp{#1}%
    \edef\sectiontitle{\expandafter\gobbleeight\meaning\temp}%
    \edef\textoutput{\write\tdm{\noindent{\noexpand\mytdmheadfont\quad\headinfo\ - \sectiontitle\hfill\noexpand\folio}\par}}%
    \textoutput%
\fi%
\font\headfontA=cmbx10 at 14pt%
{\headfontA\noindent\headinfo\ -\ #1.\hfil}%
\nobreak\vskip .5cm%
}%


\def\newsubhead#1%
{%
\ifhmode%
    \mypar%
\fi%
\ifnum\subheadno=0%
\else%
    \penalty\headpenalty\vskip .4cm%
\fi%
\nextsubheadno%
\ifmaketdm%
    \def\temp{#1}%
    \edef\subsectiontitle{\expandafter\gobbleeight\meaning\temp}%
    \edef\textoutput{\write\tdm{\noindent{\noexpand\mytdmsubheadfont\quad\quad\subheadinfo\ - \subsectiontitle\hfill\noexpand\folio}\par}}%
    \textoutput%
\fi%
\font\subheadfontA=cmsl10 at 12pt
{\subheadfontA\noindent\subheadinfo\ #1.\hfil}%
\nobreak\vskip .25cm%
}%

%
%


\font\mathromanten=cmr10
\font\mathromanseven=cmr7
\font\mathromanfive=cmr5
\newfam\mathromanfam
\textfont\mathromanfam=\mathromanten
\scriptfont\mathromanfam=\mathromanseven
\scriptscriptfont\mathromanfam=\mathromanfive
\def\mathroman{\fam\mathromanfam}


\font\sansseriften=cmss10
\font\sansserifseven=cmss10 at 7pt
\font\sansseriffive=cmss10 at 5pt
\newfam\sansseriffam
\textfont\sansseriffam=\sansseriften
\scriptfont\sansseriffam=\sansserifseven
\scriptscriptfont\sansseriffam=\sansseriffive
\def\mathsf{\fam\sansseriffam}


\font\boldten=cmb10
\font\boldseven=cmb10 at 7pt
\font\boldfive=cmb10 at 5pt
\newfam\mathboldfam
\textfont\mathboldfam=\boldten
\scriptfont\mathboldfam=\boldseven
\scriptscriptfont\mathboldfam=\boldfive
\def\mathbf{\fam\mathboldfam}


\font\mycmmiten=cmmi10
\font\mycmmiseven=cmmi10 at 7pt
\font\mycmmifive=cmmi10 at 5pt
\newfam\mycmmifam
\textfont\mycmmifam=\mycmmiten
\scriptfont\mycmmifam=\mycmmiseven
\scriptscriptfont\mycmmifam=\mycmmifive

\def\hexa#1{\ifcase #1 0\or 1\or 2\or 3\or 4\or 5\or 6\or 7\or 8\or 9\or A\or B\or C\or D\or E\or F\fi}
\mathchardef\mathi="7\hexa\mycmmifam7B
\mathchardef\mathj="7\hexa\mycmmifam7C


\font\mymsbmten=msbm10 at 8pt
\font\mymsbmseven=msbm7 at 5.6pt
\font\mymsbmfive=msbm5 at 4pt
\newfam\mymsbmfam
\textfont\mymsbmfam=\mymsbmten
\scriptfont\mymsbmfam=\mymsbmseven
\scriptscriptfont\mymsbmfam=\mymsbmfive

\mathchardef\mybeth="7\hexa\mymsbmfam69
\mathchardef\mygimmel="7\hexa\mymsbmfam6A
\mathchardef\mydaleth="7\hexa\mymsbmfam6B


\def\placelabel[#1][#2]#3{{%
\setbox10=\hbox{\raise #2cm \hbox{\hskip #1cm #3}}%
\ht10=0pt%
\dp10=0pt%
\wd10=0pt%
\box10}}%


\newif\ifinproclaim%
\global\inproclaimfalse%
\def\proclaim#1{%
\medskip%
%
%
\bgroup%
\inproclaimtrue%
\setbox10=\vbox\bgroup\leftskip=0.8em\noindent{\bf #1}\sl%
}

\def\endproclaim{%
\egroup%
\setbox11=\vtop{\noindent\vrule height \ht10 depth \dp10 width 0.1em}%
\wd11=0pt%
\setbox12=\hbox{\copy11\kern 0.3em\copy11\kern 0.3em}%
\wd12=0pt%
\setbox13=\hbox{\noindent\box12\box10}%
\noindent\unhbox13%
\egroup%
\medskip\ignorespaces%
}

\def\proclaim#1{%
\medskip%
\bgroup%
\inproclaimtrue%
\noindent{\bf #1}%
\nobreak\medskip%
\sl%
}

\def\endproclaim{%
\mypar\egroup\penalty\proclaimpenalty\medskip\ignorespaces%
}

\def\noskipproclaim#1{%
\medskip%
\bgroup%
\inproclaimtrue%
\noindent{\bf #1}\nobreak\sl%
}

\def\endnoskipproclaim{%
\mypar\egroup\penalty\proclaimpenalty\medskip\ignorespaces%
}


\def\proof{{\noindent\bf Proof:\ }}

\def\mlim{\mathop{{\mathroman Lim}}}

\def\minf{\mathop{{\mathroman Inf}}}

\def\qed{~$\square$}
\def\munion{\mathop{\cup}}
\def\minter{\mathop{\cap}}
\def\myitem#1{%
\ifinproclaim%
    \item{#1}%
\else%
    \noindent\hbox to .5cm{\hfill#1\hss}
\fi}%

\catcode`\@=11
\def\Eqalign#1{\null\,\vcenter{\openup\jot\m@th\ialign{%
\strut\hfil$\displaystyle{##}$&$\displaystyle{{}##}$\hfil%
&&\quad\strut\hfil$\displaystyle{##}$&$\displaystyle{{}##}$%
\hfil\crcr #1\crcr}}\,}
\catcode`\@=12

\def\makeop#1{%
\global\expandafter\def\csname op#1\endcsname{{\mathroman #1}}}%

\def\makeopsmall#1{%
\global\expandafter\def\csname op#1\endcsname{{\mathroman{\lowercase{#1}}}}}%

\makeopsmall{ArcTan}%
\makeopsmall{ArcCos}%
\makeop{Arg}%
\makeop{Det}%
\makeopsmall{Log}%
\makeop{Re}%
\makeop{Im}%
\makeop{Dim}%
\makeopsmall{Tan}%
\makeop{Ker}%
\makeopsmall{Cos}%
\makeopsmall{Sin}%
\makeop{Exp}%
\makeopsmall{Tanh}%
\makeop{Tr}%
\makeop{End}%
\makeop{Long}%
\makeop{Ch}%
\makeop{Exp}%
\makeop{Int}%
\makeop{Ext}%
\makeop{Aire}%
\makeop{Im}%
\makeop{Conf}%
\makeop{Mod}%
\makeop{Log}%
\makeop{Ext}%
\makeop{Int}%
\makeop{Dist}%
\makeop{Aut}%
\makeop{Id}%
\makeop{SO}%
\makeop{Homeo}%
\makeop{Vol}%
\makeop{Ric}%
\makeop{Hess}%
\makeop{Euc}%
\makeop{Isom}%
\makeop{Max}%
\makeop{Long}%
\makeop{Fixe}%
\makeop{Wind}%
\makeop{Mush}%
\makeop{Ad}%
\makeop{loc}%
\makeop{Len}%
\makeop{Area}%
\makeop{SL}%
\makeop{GL}%
\makeop{dVol}%
\makeop{Min}%
\makeop{Symm}%
\makeopsmall{Cosh}
\makeopsmall{Sinh}
\makeopsmall{ArcSinh}
\makeopsmall{ArcCotanh}
\makeopsmall{Coth}
\makeop{Ind}
\makeop{Supp}
\makeop{O}%

\let\emph=\bf

\hyphenation{quasi-con-formal}

%
%

\ifmakebiblio%
    \openout\biblio=biblio.tex %
    {%
        \edef\fileoutput{\write\biblio{\bgroup\leftskip=2em}}%
        \fileoutput
    }%
\fi%

\newref{BallGromSch}{Ballman W., Gromov M., Schroeder V., Manifolds of nonpositive curvature}
\newref{Grom}{Gromov M., Foliated plateau problem, part I : Minimal varieties}
\newref{LabC}{Labourie F., Probl\`eme de Minkowski et surfaces \`a courbure constante dans les 
vari\'et\'es hyperboliques}
\newref{LabB}{Labourie F., Probl\`emes de Monge-Amp\`ere, courbes holomorphes et laminations}
\newref{LabA}{Labourie F., Un lemme de Morse pour les surfaces convexes}
\newref{RosSpruck}{Rosenberg H., Spruck J. On the existence of convex hyperspheres of constant Gauss curvature in hyperbolic space}
\newref{SmiB}{Smith G., Hyperbolic Plateau Problems}
\newref{SmiC}{Smith G., Pointed k-surfaces}
\newref{SmiE}{Smith G., Th\`ese de doctorat, Paris (2004)}%

\ifmakebiblio%
    {\edef\fileoutput{\write\biblio{\egroup}}%
    \fileoutput}%
\fi%

%
%
%
\document
\myfontdefault
\global\chapno=1
\global\showpagenumflag=1
\def\Pagetitle{}
\null
\vfill
\def\centre{\rightskip=0pt plus 1fil \leftskip=0pt plus 1fil \spaceskip=.3333em \xspaceskip=.5em \parfillskip=0em \parindent=0em}%
\def\textmonth#1{\ifcase#1\or January\or February\or March\or April\or May\or June\or July\or August\or September\or October\or November\or December\fi}
\font\abstracttitlefont=cmr10 at 14pt
{\abstracttitlefont\centre Finite Area and Volume of Pointed $k$-Surfaces\par}
\bigskip
{\centre Graham Smith\par}
\bigskip
{\centre \the\day\ \textmonth\month\ \the\year\par}
\bigskip
{\centre Max Planck Institute for Mathematics in the Sciences,\par
Inselstrasse 22.,\par
D-04103 Leipzig,\par
GERMANY\par}
\bigskip
\noindent{\emph Abstract:\ }We define the ``volume'' contained by pointed $k$-surfaces, first studied by the 
author in \cite{SmiE}, and we show that this volume is always finite. Likewise, we show that the surface area of a pointed $k$-surface is always finite.
\vfill
\nextoddpage
\def\Pagetitle{\sl Finite Area and Volume of Pointed $k$-Surfaces}
\newhead{Introduction}
\noindent Immersed hypersurfaces of constant Gaussian curvature are very classical objects of study which in
recent years (in geometric terms) have found various fruitful applications to the study of negatively curved 
manifolds. In \cite{RosSpruck}, Rosenberg and Spruck constructed a large class of examples by solving the 
boundary value problem. Labourie then showed in \cite{LabB} that when the ambient manifold is three dimensional, constant Gaussian curvature surfaces may be studied in terms of pseudo-holomorphic curves
in a contact manifold. The powerful techniques of this latter theory then allow the construction 
\cite{LabA} of a much more general family of such surfaces, which are well adapted to a number of
useful applications, such as the construction \cite{LabC} of a canonical foliation of the non-compact ends of certain hyperbolic manifolds; the realisation of homomorphisms of compact Fuchsian groups into Kleinian
groups as constant Gauss curvature immersions (ch.4 of \cite{SmiE}); and the canonical association \cite{SmiC} of a complete immersed surface in $\Bbb{H}^3$ to each ramified covering of the Riemann sphere, which is the case that we study in this paper.
\medskip
\noindent Let $\Bbb{H}^3$ be three dimensional hyperbolic space. We identify the ideal boundary of $\Bbb{H}^3$ with
the Riemann sphere $\hat{\Bbb{C}}$. Let $\Sigma$ be a compact Riemann surface. Let $\Cal{P}$ be a finite subset of 
$\Sigma$ and denote $\Sigma'=\Sigma\setminus\Cal{P}$. Let $\varphi:\Sigma\rightarrow\hat{\Bbb{C}}$ be a ramified covering
with critical points contained in $\Cal{P}$. The pair $(\Sigma',\varphi)$ defines a Plateau problem in the sense of
\cite{LabA}. Since $\Sigma'$ is of hyperbolic type, by \cite{SmiB}, for all $k\in]0,1[$, there exists a unique immersion
$i_k:\Sigma'\rightarrow\Bbb{H}^3$ of constant Gaussian curvature equal to $k$ which is a solution to this Plateau problem
(see section \headref{Asymptotically_Tubular_Immersed_Surfaces}). In \cite{SmiC} we completely described the geometry of 
the immersed surface $(\Sigma',i_k)$, showing that it is complete and asymptotically tubular of finite order 
(in the sense of \cite{SmiC}) near the critical points (see section \headref{Asymptotically_Tubular_Immersed_Surfaces}). Heuristically, the immersed surface has only a finite number of point singularites, all of which wrap a finite number of times in a cusp shaped manner about a geodesic.
\medskip
\noindent The aim of this paper is to study the area of and the ``volume'' contained by the immersed surface 
$(\Sigma',i_k)$. The area is a relatively trivial matter, and we obtain:
\proclaim{Theorem \nextprocno}
\noindent The area of $(\Sigma',i_k)$ is finite.
\proclabel{Area_Is_Finite}
\endproclaim
\noindent We volume is more subtle, since $(\Sigma',i_k)$ is not embedded, and therefore does not have a well defined 
interior. Nonetheless, we may define $\opVol(\Sigma',i_k)$ by integrating primitives of the volume form over
$(\Sigma',i_k)$. This still poses difficulties, since the primitive of the volume form is not necessarily $L^1$ 
over $(\Sigma',i_k)$. However, we show that $\opVol(\Sigma',i_k)$ may be obtained as the limit of finite integrals.
Indeed, using notation from section \headref{Asymptotically_Tubular_Immersed_Surfaces}, for all $p\in\Cal{P}$, let 
$\gamma_p$ be a central geodesic for $(\Sigma',i_k)$ at $p$. Let $(\alpha_p,\Omega_p)$ be an asymptotically tubular chart 
for $(\Sigma',i_k)$ about $\gamma_p$ and let $f_p$ be the graph function of $(\Sigma',i_k)$ over this chart. For all 
$p\in\Cal{P}$ and for all $t>0$, we define $c_{t,p}(s) = f_p(t,s)$. For all $t>0$, we define $\Sigma'_t$ by:
$$
\Sigma'_t = \Sigma'\setminus\munion_{p\in\Cal{P}}\alpha_p(S^1\times ]t,+\infty[).
$$
\noindent Let $\beta$ be any primitive of the volume form of $\Bbb{H}^3$. Let 
$\Psi_{\gamma_p}:\Bbb{R}\times\Bbb{R}^2\rightarrow\Bbb{H}^3$ be the parameterisation given by polar
coordinates about the geodesic $\gamma_p$, is in section \headref{Polar_Coordinates_About_A_Geodesic}.
We define $\opVol(\Sigma',i_k;t)$ by:
$$
\opVol(\Sigma',i_k;t) = \int_{\Sigma'_t}i_k^*\beta
+ \sum_{p\in\Cal{P}}\int_{\{t\}\times\Bbb{R}^2}\opWind(c_{t,p},x)\Psi_{\gamma_p}^*\beta.
$$
\noindent We obtain the following result:
\proclaim{Theorem \nextprocno}
\noindent $\opVol(\Sigma',i_k;t)$ converges to a finite limit as $t$ tends to $+\infty$. Moreover, this limit is 
independant of the choices of $(\gamma_p)_{p\in\Cal{P}}$, $(\alpha_p,\Omega_p)_{p\in\Cal{P}}$ or $\alpha$.
\proclabel{Volume_Is_Finite}
\endproclaim
\noindent This now allows us to define the volume of $(\Sigma',i_k)$ as follows:
$$
\opVol(\Sigma',i_k) = \mlim_{t\rightarrow+\infty}\opVol(\Sigma',i_k;t).
$$
\noindent The key step in proving both these results lies in showing that the cusp ends of $(\Sigma',i_k)$
taper off exponentially fast. Following the same philosophy as in \cite{SmiC}, this is acheived by proving the result for the case $(D\setminus\{p\},i)$ of a complete immersed disc with a unique singularity in its
interior, and then showing that the general result may be deduced from this case.
\medskip
\noindent These results allow us to define two interesting functions over the Teichm\"uller space of ramified coverings
over the sphere, and provoke the following natural questions:
\medskip
\myitem{(i)} How do the volume and the area vary as a function of $\varphi$?
\medskip
\myitem{(ii)} What is the asymptotic behaviour of the volume and the area as $k$ tends to $0$ or $1$?
\medskip
\myitem{(iii)} How do the higher coefficents of the asymptotic series of the volume and the area vary as a function of $\varphi$?
\medskip
\myitem{(iv)} How may these new functions be related to other known functions over the Teichm\"uller space?
\medskip
\noindent This paper is structured as follows: in sections \headref{Polar_Coordinates_About_A_Geodesic}
and \headref{Surfaces_Of_Revolution} we calculate polar coordinates about a geodesic in $\Bbb{H}^3$ in order to determine the Gaussian curvature of an arbitrary surface of revolution about that geodesic. In
section \headref{Surfaces_Of_Constant_Curvature}, we use the resulting ODE to determine the asymptotic
behaviour of a constant Gaussian curvature surface of revolution about a geodesic of some function $f$. Although we only need to know the decay rates of $f$ and $f'$ in order to control the volume and the area
respectively, for no extra effort, we are also able to determine the decay rates of every derivative of $f$.
In section \headref{Asymptotically_Tubular_Immersed_Surfaces}, we recall the notion of immersed surfaces
being tubular near a critical point, as defined in \cite{SmiC}, and we adapt this notion to the current 
context. Finally, using elementary properties of convex curves in $\Bbb{R}^2$ obtained in section
\headref{Convex_Curves_In_Real_And_Hyperbolic_Space}, we prove in section 
\headref{The_Volume_Contained_By_A_Cusp} the finiteness of the volume integral over each cusp, which 
allows us in section \headref{Finiteness_Of_Area_And_Volume} to rapidly deduce Theorems
\procref{Area_Is_Finite} and \procref{Volume_Is_Finite}.
\medskip
\noindent I am grateful to Fran\c{c}ois Labourie for introducing me to the study of Plateau problems and to Jean-Marc
Schlenker for encouraging me to address this aspect of their geometry.
\newhead{Polar Coordinates About a Geodesic}
\noindent Let $\Bbb{H}^3$ be three dimensional hyperbolic space and let $g=g_{ij}$ be the hyperbolic metric. We begin by 
calculating $g$ in terms of polar coordinates about a geodisic. We identify $\Bbb{H}^3$ with the three dimensional upper 
half space:
\headlabel{Polar_Coordinates_About_A_Geodesic}
$$\matrix
\Bbb{H}^3 \hfill&= \left\{(x,y,t)\in\Bbb{R}^3\text{ s.t. }t>0\right\},\hfill\cr
g_{ij} \hfill&= t^{-2}\delta_{ij}.\hfill\cr
\endmatrix$$
\noindent Let $\gamma:\Bbb{R}\rightarrow\Bbb{H}^3$ be a geodesic. By applying an isometry of $\Bbb{H}^3$, we may assume 
that $\gamma$ is the unique geodesic going from $0$ to $\infty$ and that $\gamma(0)=(0,0,1)$. Let $N\gamma$ be the normal 
bundle over $\gamma$. We identify $N\gamma_0$, the fibre over $0$, isometrically with $\Bbb{R}^2$. Using parallel transport,
we obtain a bundle isometry $\tau_\gamma:\Bbb{R}\times\Bbb{R}^2\rightarrow N\gamma$. Let $\opExp:T\Bbb{H}^3\rightarrow\Bbb{H}^3$ be the 
exponential map. We now define $\Phi_\gamma$ by:
$$
\Phi_\gamma = \opExp\circ\tau_\gamma.
$$
\noindent This mapping is unique up to translation of the $\Bbb{R}$ coordinate and rotation of the
$\Bbb{R}^2$ coordinate. Using polar coordinates of the $\Bbb{R}^2$ component, $\Phi_\gamma$ is given explicitely by:
$$
\Phi_\gamma(t,r,\theta)=(e^t\opTanh(r)\opCos(\theta),e^t\opTanh(r)\opSin(\theta),e^t\opCosh(r)^{-1}).
$$
\noindent We have the following result:
\proclaim{Lemma \nextprocno}
\noindent With respect to the basis $(\partial_t,\partial_r,\partial_\theta)$, the metric $\Phi_\gamma^* g$ is given by:
$$
\Phi_\gamma^*g = 
\pmatrix
\opCosh^2(r)\hfill & & \cr
 & 1\hfill& \cr
 & & \opSinh^2(r)\hfill \cr
\endpmatrix.
$$
\endproclaim
\proof The vectors $\partial_t$ point along lines defined by $r$ and $\theta$ being constant. These are straight lines in
$\Bbb{R}^3$ leaving the origin. Likewise, the vectors $\partial_r$ point along vertical circles in $\Bbb{R}^3$ having the 
origin as their centre. Finally, the vectors $\partial_\theta$ point along horizontal circles in $\Bbb{R}^3$ having their 
origin on the vertical line $t\mapsto(0,0,t)$ which passes through the origin. These vectors are pairwise 
orthogonal in $\Bbb{R}^3$. Since the hyperbolic metric of $\Bbb{H}^3$ is conformally equivalent to the 
hyperbolic metric of $\Bbb{R}^3$, it follows that these vectors are also orthogonal in $\Bbb{H}^3$. It now remains to calculate the lengths 
of these vectors with respect to the hyperbolic metric.
\medskip
\noindent We calculate these vectors over the point 
$(e^t\opTanh(r)\opCos(\theta),e^t\opTanh(r)\opSin(\theta),e^t\opCosh(r)^{-1})$. Firstly:
$$
\partial_t = (e^t\opTanh(r)\opCos(\theta),e^t\opTanh(r)\opSin(\theta),e^t\opCosh(r)^{-1}).
$$
\noindent Thus:
$$\matrix
\|\partial_t\|^2 \hfill&
= e^{-2t}\opCosh^2(r)(e^{2t}\opTanh^2(r)\opCos^2(\theta) + 
e^{2t}\opTanh^2(r)\opSin^2(\theta) +
e^{2t}\opCosh^{-2}(r))\hfill\cr
&=\opCosh^2(r).\hfill\cr
\endmatrix
$$
\noindent The result now follows by an analogous calculation for $\partial_r$ and $\partial_\theta$.\qed
\medskip
\noindent We define the mapping $\Psi_\gamma:\Bbb{R}\times\Bbb{R}^+\times[0,2\pi]\rightarrow\Bbb{H}^3$ by:
$$
\Psi_\gamma(t,R,\theta) = \Phi(t,\opArcSinh(R),\theta).
$$
\noindent This mapping also yields a form of polar coordinates for $\Bbb{H}^3$ about a geodesic. However, the corresponding
metric has a simpler formula, as the following lemma shows:
\proclaim{Lemma \nextprocno}
\noindent With respect to the basis $(\partial_t,\partial_R,\partial_\theta)$, the metric $\Psi_\gamma^* g$ is given by:
\proclabel{Metric_In_Polar_Coordinates}
$$
\Phi_\gamma^*g = 
\pmatrix
(1+R^2)\hfill & & \cr
 & (1+R^2)^{-1}\hfill& \cr
 & & R^2\hfill \cr
\endpmatrix.
$$
\endproclaim
\proof Since $\partial_R$ is merely a rescaling of $\partial_r$, the pairwise orthogonality of the three coordinate vectors
is preserved by this reparametrisation. Moreover, the lengths of the vectors $\partial_t$ and $\partial_\theta$ also 
remain unchanged. Finally, since $r=\opSinh(R)$, we have:
$$\matrix
&\partial_r \hfill&= \opCosh(R)\partial_R\hfill\cr
\Rightarrow\hfill &\|\partial_R\|^2 \hfill&= \opCosh(R)^{-2}\|\partial_r\|^2.\hfill\cr
\endmatrix
$$
\noindent Using the classical relation $\opCosh^2(x) - \opSinh^2(x)=1$, we obtain the desired result.\qed
\medskip
\noindent We now calculate the action on this basis of the Levi-Civita covariant derivative of $\Psi_\gamma^*g$. We obtain 
the following result:
\proclaim{Lemma \nextprocno}
\noindent The Levi-Civita covariant derivative of $\Psi_\gamma^*g$ is determined by the following relations:
\proclabel{Connection_One_Form_In_Polar_Coordinates}
$$
\matrix
\nabla_{\partial_t}\partial_t \hfill&= -R(1+R^2)\partial_R,\hfill\cr
\nabla_{\partial_t}\partial_R \hfill&= R(1+R^2)^{-1}\partial_t,\hfill\cr
\nabla_{\partial_t}\partial_\theta \hfill&= 0,\hfill\cr
\nabla_{\partial_R}\partial_R \hfill&= -R(1+R^2)^{-1}\partial_R,\hfill\cr
\nabla_{\partial_R}\partial_\theta \hfill&= R^{-1}\partial_\theta,\hfill\cr
\nabla_{\partial_\theta}\partial_\theta \hfill&= -R(1+R^2)\partial_R.\hfill\cr
\endmatrix
$$
\endproclaim
\proof This follows directly from the preceeding lemma and the Kozhul formula.\qed
\newhead{Surfaces of Revolution}
\noindent Let $I$ be an interval in $\Bbb{R}$. Let $f:I\rightarrow]0,\infty[$ be a positive valued smooth function. 
We define $\Sigma_{f,\gamma}\subseteq\Bbb{H}^3$ by:
\headlabel{Surfaces_Of_Revolution}
$$
\Sigma_{f,\gamma} = \left\{ \Psi_\gamma(t,f(t),\theta)\text{ s.t. }t\in I,\theta\in [0,2\pi]\right\}.
$$
\noindent $\Sigma_{f,\gamma}$ is a surface of revolution in $\Bbb{H}^3$ about the geodesic $\gamma$. We aim to 
obtain differential conditions on $f$ for the surface $\Sigma_f$ to have constant Gaussian curvature. Let $\kappa(t)$
be the Gaussian curvature of the surface $\Sigma_f$ at the point $\Psi(t,f(t),0)$. We have the following result:
\proclaim{Lemma \nextprocno}
\noindent The Gaussian curvature, $\kappa$ satisfies:
\proclabel{Gaussian_Curvature_Of_Surface_Of_Revolution}
$$
\kappa f((1+f^2)+(f')^2(1+f^2)^{-1})^{3/2} = (1+f^2)(-f''(1+f^2) + f(1+f^2)^2 + 3f(f')^2). 
$$
\endproclaim
\proof We work now in the coordinates of $\Bbb{R}\times ]0,\infty[\times[0,2\pi]$. We define the function
$\hat{f}:I\times[0,2\pi]\rightarrow\Bbb{R}\times ]0,\infty[\times[0,2\pi]$ by:
$$
\hat{f}(t,\theta) = (t,f(t),\theta).
$$
\noindent We define the vector fields $\hat{\partial}_t = D\hat{f}\cdot\partial_t$ and 
$\hat{\partial}_\theta=D\hat{f}\cdot\partial_\theta$. These vector fields span the tangent space of $\Sigma_f$. We have:
$$\matrix
\hat{\partial}_t(t,\theta) \hfill&= (1,f'(t),0),\hfill\cr
\hat{\partial}_\theta(t,\theta) \hfill&= (0,0,1).\hfill\cr
\endmatrix
$$
\noindent We now define the vector field $\hat{\mathsf{N}}$ by:
$$
\hat{\mathsf{N}}(t,\theta) = (-f',(1+f^2)^2,0).
$$
\noindent This vector field spans the normal bundle to $\Sigma_f$. Moreover:
$$
\|\hat{\mathsf{N}}\|^2 = (f')^2(1+f^2) + (1+f^2)^3.
$$
\noindent By taking the covariant derivative of this vector field with respect to $\hat{\partial}_t$ and 
$\hat{\partial}_\theta$, we obtain the second fundamental form of $\Sigma_f$. Let $D$ be the canonical flat 
connexion of $\Bbb{R}\times ]0,\infty[\times[0,2\pi]$. We have:
$$\matrix
D_{\hat{\partial}_t}\hat{\mathsf{N}} \hfill&= \partial_t\hat{\mathsf{N}} \hfill&= (-f'',4(1+f^2)ff',0), \hfill\cr
D_{\hat{\partial}_\theta}\hat{\mathsf{N}} \hfill&= \partial_\theta\hat{\mathsf{N}} \hfill&= (0,0,0). \hfill\cr
\endmatrix
$$
\noindent By Lemma \procref{Metric_In_Polar_Coordinates}, we have:
$$\matrix
\langle D_{\hat{\partial}_t}\hat{\mathsf{N}}, \hat{\partial}_t \rangle \hfill&= -f''(1+f^2) + 4f(f')^2,\hfill&\qquad&
\langle D_{\hat{\partial}_t}\hat{\mathsf{N}}, \hat{\partial}_\theta \rangle \hfill&= 0,\hfill\cr
\langle D_{\hat{\partial}_\theta}\hat{\mathsf{N}}, \hat{\partial}_t \rangle \hfill&= 0,\hfill&\qquad&
\langle D_{\hat{\partial}_\theta}\hat{\mathsf{N}}, \hat{\partial}_\theta \rangle \hfill&= 0.\hfill\cr
\endmatrix$$
\noindent  Let $\Omega$ be the connexion one form of $\nabla$ with respect to $D$, so that, for any vector fields 
$X$ and $Y$:
$$
\nabla_X Y = D_X Y + \Omega(X,Y).
$$
\noindent Using Lemmata \procref{Metric_In_Polar_Coordinates} and \procref{Connection_One_Form_In_Polar_Coordinates}. We
obtain:
$$\matrix
\langle\Omega(\hat{\partial}_t,\hat{\mathsf{N}}),\hat{\partial}_t\rangle \hfill&= f(1+f^2)^2 - f(f')^2,\hfill&\qquad&
\langle\Omega(\hat{\partial}_t,\hat{\mathsf{N}}),\hat{\partial}_\theta\rangle \hfill&= 0,\hfill\cr
\langle\Omega(\hat{\partial}_\theta,\hat{\mathsf{N}}),\hat{\partial}_t\rangle \hfill&= 0,\hfill&\qquad&
\langle\Omega(\hat{\partial}_\theta,\hat{\mathsf{N}}),\hat{\partial}_\theta\rangle \hfill&= f(1+f^2)^2.\hfill\cr
\endmatrix$$
\noindent Let $\mathsf{N}$ be the unit normal vector field to $\Sigma_f$:
$$
\mathsf{N} = \|\hat{\mathsf{N}}\|^{-1}\hat{\mathsf{N}}.
$$
\noindent Let $II$ be the second fundamental form of $\Sigma_f$. That is, if $X$ and $Y$ are vector fields tangent to
$\Sigma_f$:
$$
II(X,Y) = \langle \nabla_X{\mathsf N},Y\rangle.
$$
\noindent If $X$ and $Y$ are both vector fields tangent to $\Sigma_f$, then 
$\langle\nabla_X\mathsf{N},Y\rangle=\|\hat{{\mathsf N}}\|^{-1}\langle\nabla_X\hat{{\mathsf N}},Y\rangle$. We may thus
calculate $II$:
$$\matrix
II(\hat{\partial}_t,\hat{\partial_t}) \hfill&= ((f')^2(1+f^2) + (1+f^2)^3)^{-1/2}\times\hfill\cr
&\qquad\qquad (-f''(1+f^2) + f(1+f^2)^2 + 3f(f')^2),\hfill\cr
II(\hat{\partial}_t,\hat{\partial_\theta}) \hfill&= 0,\hfill\cr
II(\hat{\partial}_\theta,\hat{\partial}_t) \hfill&= 0,\hfill\cr
II(\hat{\partial}_\theta,\hat{\partial}_\theta) \hfill&= ((f')^2(1+f^2) + (1+f^2)^3)^{-1/2}f(1+(f')^2)^2.\hfill\cr
\endmatrix$$
\noindent Observing that $\hat{\partial}_t$ and $\hat{\partial}_\theta$ are orthogonal to one another, we obtain:
$$\matrix
\opDet(\hat{\partial}_t,\hat{\partial}_\theta) \hfill&= \|\hat{\partial}_t\|^2\|\hat{\partial}_\theta\|^2 \hfill\cr
&= f^2(1+f^2)((1+f^2)+f'(1+f^2)^{-2}).\hfill\cr
\endmatrix$$
\noindent If we denote by $A$ the matrix of $II$ with respect to the basis $(\hat{\partial}_t,\hat{\partial}_\theta)$, then
the Gaussian curvature, $\kappa$, satisfies:
$$
\kappa = \opDet(A)/\opDet(\hat{\partial}_t,\hat{\partial}_\theta).
$$
\noindent Thus:
$$\matrix
\kappa f^2(1+f^2)((1+f^2)+(f')^2(1+f^2)^{-1})^{3/2} \hfill\cr
\qquad\qquad\qquad= f(1+f^2)^2(-f''(1+f^2) + f(1+f^2)^2 + 3f(f')^2).\hfill\cr
\endmatrix$$
\noindent The result now follows.\qed
\newhead{Surfaces Of Constant Curvature}
\noindent We now study the asymptotique behaviour of solutions to the differential differential equation given by 
Lemma \procref{Gaussian_Curvature_Of_Surface_Of_Revolution}. We have the following result:
\headlabel{Surfaces_Of_Constant_Curvature}
\proclaim{Lemma \nextprocno}
\noindent Let $k$ be a real number in $]0,1[$. Let $f:[0,\infty[\rightarrow]0,\infty[$ be such that the surface of
revolution $\Sigma_{f,\gamma}$ is of constant Gaussian curvature equal to $k$. Suppose, moreover, that $f(t)$ and $f'(t)$ 
both tend to zero as $t$ tends to $+\infty$. Then, for all $\delta>0$, there exists $T>0$ and constants $B>A>0$ such 
that for all $t>T$:
\proclabel{The_Function_Decays_Exponentially}
$$
Ae^{-(\lambda + \delta)t} \leqslant f(t),-f'(t),f''(t) \leqslant Be^{-(\lambda-\delta)t},
$$
\noindent where $\lambda^2 = 1-k$.
\endproclaim
\proof By Lemma \procref{Gaussian_Curvature_Of_Surface_Of_Revolution}, $f$, satisfies the following differential
equation:
$$
f''/f = (1-k) + \epsilon,
$$
\noindent where $\epsilon:[0,\infty[\rightarrow ]0,\infty[$ is a smooth function such that $\epsilon(t)$ tends to zero
as $t$ tends to infinity. We define the function $g(t) = \opLog(f(t))$. Thus:
$$
g' = f'/f,\qquad g'' = (ff'' - (f')^2)/f^2 = f''/f - (g')^2.
$$
\noindent The function $g$ therefore satisfies the following differential relation:
$$
g'' + (g')^2 = (1-k) + \epsilon.
$$
\noindent We define $h(t)=g'(t)$. We then obtain:
$$
h' + h^2 = \lambda^2 + \epsilon.
$$
\noindent Let $\delta>0$ be such that $\delta<\lambda$. Let $T_0>0$ be such that for $t>T_0$:
$$
\epsilon(t) < \delta^2.
$$
\noindent Let $t>T_0$ be arbitrary and suppose that $h(t)\geqslant \lambda+\delta$. Then:
$$
\left|\lambda^2 - h(t)^2\right| \geqslant \left|\delta^2 + 2\delta\lambda\right| \geqslant \delta\lambda. 
$$
\noindent Thus:
$$
\left|\epsilon/(\lambda^2 - h(t)^2)\right| \leqslant \delta/\lambda.
$$
\noindent Combining this with the differential equation for $h$, we obtain:
$$
\left|h'(t)(\lambda^2 - h(t)^2)^{-1} - 1\right| \leqslant \delta/\lambda.
$$
\noindent Consequently, if we define $\eta(t) = \lambda^{-1}\opArcCotanh(h(t)/\lambda)$, we obtain:
$$
\left|\eta'(t) - 1\right| \leqslant \delta/\lambda.
$$
\noindent It follows that, for any $t_1>t_0>T_0$, if $h(t)\geqslant \lambda+\delta$ for all $t$ in the 
interval $[t_0,t_1]$, then:
$$
\eta(t_1) \leqslant \eta(t_0) + (t_1 - t_0)(1+\delta/\lambda).
$$
\noindent Thus, under the same conditions:
$$
h(t_1) \leqslant \lambda\opCoth(\lambda\eta(t_0) + (\lambda+\delta)(t_1-t_0)).
$$
\noindent This tells us that if $t>T_0$ and $h(t)$ is ever greater than $\lambda+\delta$, then, in finite time, 
it will fall below $\lambda+\delta$. Moreover, for all $t>T_0$, $h'(t)<0$ whenever $h(t)=\lambda+\delta$. It thus follows
that the function will never thereafter be greater than $\lambda+\delta$. Heuristically, we have shown that the function
behaves like the positive branch of the hyperbolic cotangent. We thus conclude that, for any solution, $h$, 
there exists $T_1>0$ such that, for $t>T_1$:
$$
h(t) \leqslant \lambda+\delta.
$$
\noindent A similar analysis for the interval $]-\infty,-\lambda-\delta[$ reveals that there exists $T_2>0$ which depends
only on the function $\epsilon$ such that if there exists $t>T_2$ with $h(t)<-\lambda-\delta$, then $h(t)$ tends to $-\infty$
in finite time (since, heuristically, it behaves like the negative branch of the hyperbolic cotangent). However, $f$ 
exists and is positive for all time. Thus $g$ and $g'=h$ both exist for all time, and this is not possible. It thus 
follows that, for $t>T_2$:
$$
h(t) \geqslant -\lambda-\delta.
$$
\noindent Finally, by considering the interval $]-\lambda+\delta,\lambda-\delta[$, there exists $T_3>0$ which depends only
on $\epsilon$ and a constant $\Delta T_3>0$ which depends only on $\delta$ such that if there exists $t>T_3$ with
$h(t)>-\lambda+\delta$, then, for all $t'>t + \Delta T_3$, $h(t')>\lambda - \delta$. This happens heuristically because the
solution in this case behaves like the hyperbolic tangent.
\medskip
\noindent We have thus shown that there exists $T_4>0$ such that, for all $t>T_4$, either $\left|h(t)-\lambda\right|<\delta$
or $\left|h(t)+\lambda\right|<\delta$.
\medskip
\noindent We now exclude the case where $\left|h(t)-\lambda\right|<\delta$ when $t>T_4$. Indeed, suppose that
$\delta<\lambda/2$. In this case, since $h(t)=g'(t)$, it follows that, for large values of $t$, the function 
$g(t)$ grows faster than $\lambda t/2$. Since $f$ is the exponential of $g$ it then follows that $f$ tends to
infinity as $t$ tends to infinity, and this contradicts the hypotheses on $f$.
\medskip
\noindent It thus follows that, for all $t>T_4$, $\left|h(t) + \lambda\right|<\delta$. Consequently, there exists a constant 
$C_1$ such that, for all $t>T_4$:
$$
C_1 - (\lambda +\delta) t \leqslant g(t) \leqslant C_1 - (\lambda-\delta) t.
$$
\noindent Taking the exponential of each of these functions, we see that there exists a constant $C_2$ such that, for $t>T_4$:
$$
C_2 e^{-(\lambda + \delta)t} \leqslant f(t) \leqslant C_2e^{-(\lambda+\delta)t}.
$$
\noindent Since $f''=f((1-k)+\epsilon)$, there exist constants $C_3,C_4>0$ and $T_5\geqslant T_4$ such that for $t>T_5$:
$$
C_3 e^{-(\lambda + \delta)t} \leqslant f''(t) \leqslant C_4e^{-(\lambda+\delta)t}.
$$
\noindent Finally, since $f'(t)$ tends to $0$ as $t$ tends to $+\infty$, we obtain the relation for $f'(t)$ by integrating
$f''(t)$ back from $+\infty$. The result now follows.\qed
\medskip
\noindent We may also estimate the higher derivatives of $f$:
\proclaim{Corollary \nextprocno}
\noindent With the hypothesis and notation of the previous lemma, for all $k\geqslant 2$, there exists $B_k>A_k>0$ such 
that for $t\geqslant T$:
$$
A_k e^{-(\lambda+\delta)t} \leqslant (-1)^kf^{(k)}(t) \leqslant B_k e^{-(\lambda-\delta)}t.
$$
\endproclaim
\proof By induction, for all $k\geqslant 0$:
$$
f^{(k+2)} = f^{(k)}(1-k) + \Sigma_{i=0}^k\epsilon_if^{(i)},
$$
\noindent where, for all $i$, $\epsilon_i(t)$ tends to zero as $t$ tends to $+\infty$. The result now follows by induction.\qed
\medskip
\noindent This allows us to control the area of $\Sigma_{f,\gamma}$ and the volume that it contains:
\proclaim{Corollary \nextprocno}
\noindent Let $\opArea(t)$ and $\opVol(t)$ be respectively that area of and the volume inside the restriction of
$\Sigma_f$ to $[t,+\infty[$. Then, for all $\delta>0$, there exists $T>0$ and $B>A>0$ such that, for all $t>T$:
\proclabel{Control_Of_Area_And_Volume_Of_Rotationally_Symmetric_Cusp}
$$\matrix
Ae^{-(\lambda+\delta)t} \hfill&\leqslant \opArea(t) \hfill\leqslant Be^{-(\lambda+\delta)t},\hfill\cr
Ae^{-2(\lambda+\delta)t} \hfill&\leqslant \opVol(t) \hfill\leqslant Be^{-2(\lambda+\delta)t},\hfill\cr
\endmatrix$$
\noindent where $\lambda^2 = 1-k$. 
\endproclaim
\proof This follows directly by calculating the area and volume integrals, bearing in mind that $\Psi_\gamma^*g$ is
uniformly equivalent to the Euclidian metric in an $\epsilon$-neighbourhood of $\gamma$.\qed
\newhead{Asymptotically Tubular Immersed Surfaces}
\noindent Let $\Sigma$ be a compact Riemann surface and let $\Cal{P}$ be a finite subset of $\Sigma$. Define 
$\Sigma'=\Sigma\setminus\Cal{P}$. Let $\varphi:\Sigma\rightarrow\hat{\Bbb{C}}$ be a ramified covering such that the 
ramification points are contained in $\Cal{P}$. The pair $(\Sigma',\varphi)$ defines a Plateau problem in the sense
of Labourie, \cite{LabA}. 
\headlabel{Asymptotically_Tubular_Immersed_Surfaces}
\medskip
\noindent Let $i:\Sigma\rightarrow\Bbb{H}^3$ be a convex immersion. Let ${\mathsf N}_i:\Sigma\rightarrow U\Bbb{H}^3$ be 
the exterior unit normal over $i$. We call this the Gauss lifting of $i$ and in the sequel we denote it by
$\hat{\mathi}$. Let $\overrightarrow{n}:U\Bbb{H}^3\rightarrow\partial_\infty\Bbb{H}^3=\hat{\Bbb{C}}$ be the 
Gauss-Minkowski mapping. Thus, if $\gamma:\Bbb{R}\rightarrow\Bbb{H}^3$ is a unit speed geodesic in $\Bbb{H}^3$, then:
$$
\overrightarrow{n}(\partial_t\gamma) = \gamma(+\infty).
$$
\noindent Since $i$ is convex, elementary hyperbolic geometry (see for example, \cite{BallGromSch}) allows us to show that $\overrightarrow{n}\circ\hat{\mathi}$
is a local homeomorphism. 
\medskip
\noindent For $k\in]0,1[$, following \cite{LabA}, the pair $(\Sigma',i)$ is said to be a solution of the Plateau problem 
$(\Sigma',\varphi)$ with Gaussian curvature equal to $k$ if and only if:
\medskip
\myitem{(i)} the mapping $i$ is a convex immersion with Gaussian curvature equal to $k$,
\medskip
\myitem{(ii)} $(\Sigma',\hat{\mathi})$ is complete in the sense of immersed surfaces, and
\medskip
\myitem{(iii)} $\varphi = \overrightarrow{n}\circ\hat{\mathi}$.
\medskip
\noindent Since the surface $\Sigma'$ is of hyperbolic conformal type, by \cite{SmiB}, for all $k\in]0,1[$, there exists a 
unique solution $(\Sigma',i_k)$ of the Plateau problem $(\Sigma',\varphi)$ with Gaussian curvature equal to $k$.
\medskip
\noindent Let $p$ be an arbitrary point in $\Cal{P}$. Let $n$ be the order of ramification of $\varphi$ at $p$. In \cite{SmiC}, we defined the notion of a surface being asymptotically tubular of finite order near a 
point singularity, and we showed that the immersed surface $(\Sigma',\hat{\mathi}_k)$ is
asymptotically tubular of order $n$ near $p$. This implies that there exists:
\medskip
\myitem{(i)} a geodesic $\gamma$ such that $\gamma(+\infty)=\varphi(p)$,
\medskip
\myitem{(ii)} a smooth function $f:S^1\times ]0,+\infty[\rightarrow\Bbb{R}^2$, 
\medskip
\myitem{(iii)} a neighbourhood $\Omega$ of $p$ in $\Sigma$ containing no other point of $\Cal{P}$, and
\medskip
\myitem{(iv)} a diffeomorphism $\alpha:S^1\times ]0,+\infty[\rightarrow\Omega\setminus\{p\}$,
\medskip
\noindent such that:
\medskip
\myitem{(i)} $\Psi_\gamma(t,f(s,t))=(i\circ\alpha)(s,t)$, 
\medskip
\myitem{(ii)} $\alpha(s,t)$ tends towards $p$ as $t$ tends to $+\infty$, and
\medskip
\myitem{(iii)} $f(\cdot,t+\cdot)$ converges to $0$ in the $C^\infty_\oploc$ topology as $t$ tends to 
$+\infty$.
\medskip
\noindent Moreover:
\medskip
\myitem{(iv)} for all $t$, $f(\cdot,t)$ has index $n$ in a sense that will be made clear shortly.
\medskip
\noindent In the sequel, we refer to $\gamma$ as a central geodesic for $(\Sigma,i_k)$ at $p$, we refer to
$(\alpha,\Omega)$ as an asymtotically tubular chart for $(\Sigma,i_k)$ about $\gamma$ at $p$, and we refer to $f$ as the 
graph function of $(\Sigma,i_k)$ over this chart. 
\medskip
\noindent We have the following result:
\proclaim{Lemma \nextprocno}
\noindent Let $\gamma$ be a central geodesic for $(\Sigma,i_k)$ at $p$ and let $(\alpha,\Omega)$ be an asymptotically 
tubular chart for $(\Sigma,i_k)$ about $\gamma$ at $p$. For all $t\geqslant 0$, $(\Omega\setminus\{p\},i_k)$ is transverse 
to $\Psi_\gamma(\{t\}\times\Bbb{R}^2)$.
\proclabel{Surface_Is_Transverse}
\endproclaim
\proof Let $p:\Bbb{H}^3\rightarrow\gamma$ be the orthogonal projection. Let $f$ be the graph function of $(\Sigma,i_k)$
over $(\alpha,\Omega)$. Then:
$$\matrix
&(i_k\circ\alpha)(s,t) \hfill&= \Psi_\gamma(t,f(s,t))\hfill\cr
\Rightarrow\hfill&(p\circ i_k\circ\alpha)(s,t) \hfill&= \gamma(t).\hfill\cr
\endmatrix$$
\noindent The orthogonal projection onto the geodesic is thus surjective, and the result now follows.\qed
\medskip
\noindent For $t>0$, we define the mapping $c_t=f(\cdot,t)$, and we obtain the following corollary:
\proclaim{Corollary \nextprocno}
\noindent For $t>0$, the mapping $c_t$ is a smooth immersed curve. Moreover, for all $t$, $c_t$ is convex with respect to the hyperbolic metric $\Psi_\gamma^*g$ over $\{t\}\times\Bbb{R}^2$.
\proclabel{Intersection_Is_Smooth_Curve}
\endproclaim
\proof By transversality, $c_t$ is immersed. Since $(\Sigma',i_k)$ is convex and $\{t\}\times\Bbb{R}^2$ is
totally geodesic, $c_t$ is also convex. The result now follows.\qed
\medskip
\noindent For all $t$, we define the 
exterior unit normal ${\mathsf N}_t$ of $c_t$. We then orient $c_t$ such that ${\mathsf N}_t$ lies to its right hand side.
By composing ${\mathsf N}_t$ with the Gauss-Minkowski mapping, we obtain a continuous mapping from $S^1$ into
$\partial_\infty\Psi_\gamma(\{t\}\times\Bbb{R}^2)$, which itself is homeomorphic to $S^1$ (the orientation of
$\partial_\infty\Psi_\gamma(\{t\}\times\Bbb{R}^2)$ may be explicitely specified although it is not very 
important). We thus define $\opInd(c_t)$, the index of $c_t$, by:
$$
\opInd(c_t) = \opInd(\overrightarrow{n}\circ{\mathsf N}_t).
$$
\noindent Condition $(iv)$ may now be made explicit:
\medskip
\myitem{(iv)} for all $t$, $\opInd(c_t)=n$.
\newhead{Convex Curves in Real and Hyperbolic Space}
\noindent We now require the following elementary results concerning the geometry of convex curves:
\headlabel{Convex_Curves_In_Real_And_Hyperbolic_Space}
\proclaim{Lemma \nextprocno}
\noindent Let $M$ be either $\Bbb{R}^2$ or $\Bbb{H}^2$, so that $\partial_\infty M$ is homeomorphic to $S^1$. Let $UM$ be
the unitary bundle of $M$ and let $\overrightarrow{n}:UM\rightarrow\partial_\infty M$ be the Gauss-Minkowski mapping.
\proclabel{Winding_Number_Is_Bounded_By_Index}
\medskip
\noindent Let $c:S^1\rightarrow M$ be a smooth, closed, convex curve. Let $p\in M$ be any point in the complement of the 
image of $c$. If $\opWind(c,p)$ be the winding number of $c$ about $p$, then:
$$
0 \leqslant \opWind(c,p) \leqslant \opInd(c).
$$
\endproclaim
\proof Let ${\mathsf N}$ be the exterior unit normal to $c$. We assume that $c$ is oriented so that ${\mathsf N}$ lies to
its right hand side. By deforming $c$ by a small amount, we obtain a curve $c'$ arbitrarily close to $c$ in the 
$C^{\infty}$ topology which is convex and intersects itself transversally. In particular, if $c'$ is sufficiently close to 
$c$, then:
$$
\opInd(c) = \opInd(c'),\qquad \opWind(c,p) = \opWind(c',p).
$$
\noindent We thus assume that $c$ intersects itself transversally. In particular, $c$ only intersects itself at a
finite number of points. We may therefore decompose $c$ into a finite collection $c_1,...,c_n$ of piecewise smooth, simple, 
closed curves which are convex except possibly at the apexes, where different curves join to each other. The number of
apexes of $c_i$ equals the number of distinct components of $c$ comprising $c_i$. For each $i$ let ${\mathsf N}_i$ be the
restriction of ${\mathsf N}$ to $c_i$. The $c_i$ may be labelled by the vertices of a tree, according to how they join to
each other. The leaves are then precisely the curves with only one apex. By induction from the leaves downwards, we may
show that each ${\mathsf N}_i$ only points into one of the connected components of the complement of $c_i$ (i.e. it 
does not change sign at the apexes).
\medskip
\noindent For each $i$, let $\Omega^0_i$ and $\Omega_i^\infty$ be respectively the bounded and unbounded components of the
complement of $c_i$ in $M$. Let $\hat{\Omega}_i$ be the convex hull of $c_i$ in $M$. Let $\Gamma$ be a supporting 
geodesic of $\hat{\Omega}_i$. $\Gamma$ intersects $c_i$ non trivially. By the convexity of $c_i$, we may assume that 
that it intersects $c_i$ away from the apexes. At this point of intersection, ${\mathsf N}_i$ points into the complement 
of $\hat{\Omega}_i$. Consequently, ${\mathsf N}_i$ always points into $\Omega^\infty_i$.
\medskip
\noindent It follows that, for all $i$, ${\mathsf N}_i$ points outwards from $\Omega^0_i$. Consequently, for each $i$,
$\opInd(c_i)$ brings a contribution of $+1$ to $\opInd(c)$, and $\opWind(c_i,p)$ brings a contribution of $0$ or $+1$ to
$\opWind(c,p)$. The result now follows.\qed
\medskip
\noindent We now recall the following generalisation of Stokes theorem:
\proclaim{Lemma \nextprocno}
\noindent Let $c:S^1\rightarrow\Bbb{R}^2$ be a smooth, closed curve. If $\alpha$ is a $1$-form over $\Bbb{R}^2$, then:
\proclabel{Stokes_Theorem}
$$
\int_c \alpha = \int_{\Bbb{R}^2}\opWind(c,x)d\alpha(x).
$$
\endproclaim
\proof By deforming $c$ a small amount, we obtain a curve $c'$ arbitrarily close to $c$ in the $C^{\infty}$ topology
such that $c'$ is convex and intersects itself transversally. By choosing $c'$ sufficiently close to $c$, we may assume
that $\opWind(c,\cdot)=\opWind(c',\cdot)$ except on a set of arbitrarily small measure. We may thus assume that $c$
intersects itself transversally. As in the proof of Lemma \procref{Winding_Number_Is_Bounded_By_Index}, we may decompose
$c$ into a finite collection $c_1,... c_n$ of simple closed curves. By Stokes' theorem, the 
result holds for each $c_i$, and the general result holds by additivity.\qed
\medskip
\noindent This also allows us to obtain the derivative of the winding number of a smoothly varying family of curves as 
a distribution over $\Bbb{R}^2$:
\proclaim{Lemma \nextprocno}
\noindent Let $c_t:S^1\times]-\epsilon,\epsilon[$ be a smoothly varying family of smooth curves in $\Bbb{R}^2$. If 
$\beta$ is a $2$-form in $\Bbb{R}^2$, then:
\proclabel{Derivative_Of_Winding_Number}
$$
\partial_t\int_{\Bbb{R}^2}\opWind(c_t,x)\beta(x) = \int_{c_t}i_{\partial_t c_t}\beta.
$$
\endproclaim
\proof Let $\Cal{L}$ denote the Lie derivative. Let $\beta':\Bbb{R}^2\rightarrow\Bbb{R}$ be a compactly supported $2$-form 
such that:
$$
\int_{\Bbb{R}^2}\beta' = 0.
$$
\noindent Let $\gamma$ be a primitive of $\beta'$. By Lemma \procref{Stokes_Theorem}, for all $t$, we have:
$$
\int_{c_t}\gamma = \int_{\Bbb{R}^2}\opWind(c_t,x)\beta(x).
$$
\noindent We have:
$$\matrix
\partial_t\int_{\Bbb{R}^2}\opWind(c_t,x)\beta'(x)\hfill&=\partial_t\int_{c_t}\gamma\hfill\cr
&=\int_{c_t}\Cal{L}_{\partial_t c_t}\gamma\hfill\cr
&=\int_{c_t}(d i_{\partial_t c_t} + i_{\partial_t c_t}d)\gamma\hfill\cr
&=\int_{c_t}i_{\partial_t c_t}\beta'.\hfill\cr
\endmatrix$$
\noindent By reducing $\epsilon$ if necessary, we may construct a $2$-form, $\beta_0$ such that $\opSupp(\beta_0)$ is 
disjoint from $c_t$ for all $t$, $\beta-\beta_0$ has compact support, and:
$$
\int_{\Bbb{R}^2}\beta-\beta_0 = 0.
$$
\noindent Since the integral of $\opWind(c_t,x)\beta_0(x)$ is constant, we obtain:
$$
\partial_t\int_{\Bbb{R}^2}\opWind(c_t,x)\beta(x)=\int_{c_t}i_{\partial_t c_t}\beta.
$$
\noindent The result now follows.\qed
\newhead{The Volume Contained by a Cusp}
\noindent Let $\gamma$ be a central geodesic for $(\Sigma,i_k)$ at $p$. Let $(\alpha,\Omega)$ be an asymtotically
tubular chart for $(\Sigma,i_k)$ about $\gamma$ at $p$, and let $f$ be the graph function of $(\Sigma,i_k)$ over this
chart. We begin by controlling the image of $(\Sigma,i_k)$:
\headlabel{The_Volume_Contained_By_A_Cusp}
\proclaim{Lemma \nextprocno}
\noindent For all $\delta>0$, there exists $T>0$ and $A>0$ such that, for all $t\geqslant T$:
\proclabel{Image_Lies_Inside_Cusp}
$$
\|f(s,t)\| \leqslant Ae^{-(\lambda-\delta)t},
$$
\noindent where $\lambda^2=1-k$.
\endproclaim
\proof By applying an isometry of $\Bbb{H}^3$, we may suppose that $\gamma$ is the unique geodesic in $\Bbb{H}^3$ 
joining $0$ to $\infty$. Let $D$ be a disc in $\hat{\Bbb{C}}$ centred about the origin such that no other point
in $\varphi(\Cal{P})$ lies in $D$. We may assume that $D$ has unit radius. Let $j:D\setminus\{0\}\rightarrow\Bbb{H}^3$ 
be the unique solution of the Plateau problem $(D\setminus\{q\},\opId)$ with constant Gaussian curvature equal to $k$. 
We observe that this mapping is an embedding. 
\medskip
\noindent We will show that the immersed surface $(\Sigma',i_k)$ lies entirely within the interior of 
$(D\setminus\{0\},j)$. Indeed, for $t\in ]0,1]$ we define $D_t\subseteq\Bbb{C}$ and $k_t\in ]0,1[$ by:
$$\matrix
D_t \hfill&= \{z\in\Bbb{C}\text{ s.t. } (1-t)/2<\left|z\right|<(1+t)/2\},\hfill\cr
k_t \hfill&= (1-t) + tk.\hfill\cr
\endmatrix$$
\noindent For all $t$, let $j_t:D_t\rightarrow\Bbb{H}^3$ be the unique solution to the Plateau problem $(D_t,\opId)$ with
constant Gaussian curvature equal to $k_t$. We see that $(D_t,j_t)_{t\in]0,1[}$ defines a foliation of the exterior of
$(D,j)$. There exists $\epsilon>0$ such that for $t<\epsilon$:
$$
(\Sigma,i)\minter(D_t,j_t) = \emptyset.
$$
\noindent Let us define $t_0\in ]0,1]$ by:
$$
t_0 = \minf\{t\in ]0,1]\text{ s.t. }(\Sigma,i)\minter(D_t,j_t) \neq \emptyset.\}.
$$
\noindent Suppose that $t_0<1$. By compactness, there is some point in the closure of $(\Sigma',i)$ in 
$\Bbb{H}^3\munion\hat{\Bbb{C}}$ which lies in the image of $(D_{t_0}, j_{t_0})$. Since $(D_{t_0},j_{t_0})$ does not 
intersect $\varphi(\Cal{P})$, it follows that $(\Sigma',i)$ intersects the image of $(D_{t_0}, j_{t_0})$ at some 
finite point of $\Bbb{H}^3$. However, since $(D_t,j_t)_{t\in]0,t_0[}$ forms a foliation of the exterior of 
$(D_{t_0},j_{t_0})$, $(\Sigma',i)$ lies in the interior of $(D_{t_0},j_{t_0})$. However, this is 
impossible by the geometric maximum principal (see, for example, \cite{LabA}), since the Gaussian curvature of 
$(D_{t_0},j_{t_0})$ is greater than that of $(\Sigma',i)$. Thus $t_0=1$, and $(\Sigma',i_k)$ lies in the interior of
$(D\setminus\{0\},j)$.
\medskip
\noindent The result now follows by Lemma \procref{The_Function_Decays_Exponentially}, since, by uniqueness, 
$(D\setminus\{0\},j)$ is a surface of revolution about $\gamma$.\qed
\medskip
\noindent We now obtain estimates concerning the ``volume'' bounded by $(\Omega,i)$. 
Let $\alpha$ be any primative of the volume of $\Bbb{H}^3$. We define the function 
$V:[T,+\infty[\rightarrow\Bbb{R}$ by:
$$\matrix
V(t) \hfill&= \int_{[T,t]\times S^1}(i\circ f)^*\alpha \hfill\cr
&\qquad +\int_{\{T\}\times\Bbb{R}^2}\opWind(c_T,x)\alpha(T,x) \hfill\cr
&\qquad -\int_{\{t\}\times\Bbb{R}^2}\opWind(c_t,x)\alpha(t,x).\hfill\cr
\endmatrix$$
\noindent First, we have:
\proclaim{Lemma \nextprocno}
\noindent Let $\opdVol$ be the hyperbolic volume element of $\Bbb{R}\times\Bbb{R}^2$. The function $V(t)$ satisfies:
\proclabel{V_Is_Volume}
$$
V(t) = \int_{[T,t]\times\Bbb{R}^2}\opWind(c_t,x)\opdVol(t,x).
$$
\endproclaim
\proof Let $\Cal{L}$ denote the Lie derivative. Let $\partial_t$ denote the derivative in the direction of the first 
coordinate in $\Bbb{R}\times\Bbb{R}^2$, and let $\partial_t c_t$ denote the infinitesimal variation of $c_t$. We 
recall that:
$$
\Cal{L}_{\partial_t} \alpha = d i_{\partial_t} \alpha + i_{\partial_t} d \alpha.
$$
\noindent Thus:
$$\matrix
\int_{[T,t]}\int_{\{s\}\times\Bbb{R}^2}\opWind(c_t,x)(i_{\partial_t}d\alpha)(x)ds\hfill&=
\int_{[T,t]}\int_{\{s\}\times\Bbb{R}^2}\opWind(c_t,x)(\Cal{L}_{\partial_t}\alpha)(x)ds\hfill\cr
&\qquad -\int_{[T,t]}\int_{\{s\}\times\Bbb{R}^2}\opWind(c_t,x)(di_{\partial_t}\alpha)(x)ds.\hfill\cr
\endmatrix$$
\noindent Since $d\alpha = \opdVol$:
$$
\int_{[T,t]}\int_{\{s\}\times\Bbb{R}^2}\opWind(c_t,x)(i_{\partial_t}d\alpha)(x)ds=
\int_{[T,t]\times\Bbb{R}^2}\opWind(c_t,x)\opdVol(t,x).
$$
\noindent Next, using Lemmata \procref{Stokes_Theorem} and \procref{Derivative_Of_Winding_Number}, and taking care with orientations:
$$\matrix
\int_{[T,t]}\int_{\{s\}\times\Bbb{R}^2}\opWind(c_t,x)(di_{\partial_t}\alpha)(x)ds\hfill&=
\int_{[T,t]}\int_{\{s\}\times S^1}(i\circ f)^*(i_{\partial_t}\alpha)(\theta)ds\hfill\cr
&=-\int_{[T,t]\times S^1}(i\circ f)^*\alpha \hfill\cr
&\qquad -\int_{[T,t]}\int_{\{s\}\times S^1}(i\circ f)^*(i_{\partial_t c_t}\alpha)(\theta)ds\hfill\cr
&=-\int_{[T,t]\times S^1}(i\circ f)^*\alpha \hfill\cr
&\qquad -\int_{[T,t]}\partial_t\int_{\{s\}\times\Bbb{R}^2}\opWind(c_t,x)\alpha(x)ds\hfill\cr
&\qquad +\int_{[T,t]}\int_{\{s\}\times\Bbb{R}^2}\opWind(c_t,x)(\Cal{L}_{\partial_t}\alpha)(x)ds.\hfill\cr
\endmatrix$$
\noindent Combining these relations, we obtain:
$$\matrix
\int_{[T,t]\times\Bbb{R}^2}\opWind(c_t,x)\opdVol(t,x)&=
\int_{[T,t]\times S^1}(i\circ f)^*\alpha\hfill\cr
&\qquad +\int_{[T,t]}\partial_t\int_{\{s\}\times\Bbb{R}^2}\opWind(c_t,x)\alpha(x)ds.\hfill\cr
\endmatrix$$
\noindent The result now follows by integrating the last integral.\qed
\medskip
\noindent This allows us to prove the convergence of $V(t)$:
\proclaim{Lemma \nextprocno}
\noindent The function $V(t)$ converges to a finite limit as $t$ tends to $+\infty$.
\proclabel{The_Volume_Function_Is_Bounded}
\endproclaim
\proof By the convexity of $(\Sigma',i_k)$, corollary \procref{Intersection_Is_Smooth_Curve} and Lemmata
\procref{Winding_Number_Is_Bounded_By_Index} and \procref{V_Is_Volume}, the function $V(t)$ is positive and increasing.
We recall that the hyperbolic metric $\Psi_\gamma^*g$ is uniformly equivalent to the Euclidean metric in an
$\epsilon$-neighbourhood of $\gamma$. It thus follows by Lemmata \procref{V_Is_Volume} and 
\procref{Image_Lies_Inside_Cusp} that $V$ is bounded from above. The result now follows.\qed
\newhead{Finiteness of Area and Volume}
\noindent We are now in a position to prove Theorem \procref{Volume_Is_Finite}:
\medskip
{\noindent\bf Proof of Theorem \procref{Volume_Is_Finite}:} The existence and finiteness of this limit follows from 
Lemma \procref{The_Volume_Function_Is_Bounded}. If $\alpha'$ is another primitive of the volume form, then 
$d(\alpha'-\alpha)=0$. Thus, since the homology of $\Bbb{H}^3$ vanishes in dimension higher than zero, the integral of 
$\alpha'-\alpha$ over any closed surface vanishes and this limit does not depend on the choice of $\alpha$. 
For any $p\in\Cal{P}$, two different asymptotically tubular charts about $\gamma_p$ differ only by a rotation of the $S^1$ 
coordinate and a translation of the $\Bbb{R}$ coordinate. It thus follows by Lemma \procref{The_Volume_Function_Is_Bounded} 
that this limit does not depend on the asymtotically tubular chart chosen. An analogous reasoning shows that the integral 
does not depend on the choice of the central geodesics. The result now follows.\qed
\headlabel{Finiteness_Of_Area_And_Volume}
\medskip
\noindent The finiteness of the area of $(\Sigma',i_k)$ is significantly simpler to prove:
\medskip
{\noindent\bf Proof of Theorem \procref{Area_Is_Finite}:} Since there are only a finite number of cusps, it suffices 
to prove that the area of each cusp is finite. Let $p$ be a point in $\Cal{P}$. Let $n$ be the order of ramification of 
the function $\varphi$ at $p$. We define $q=\varphi(p)$. Let $D$ be a disc in $\hat{\Bbb{C}}$ about $q$ which contains no 
other point of $\Cal{P}$. By applying an isometry of $\Bbb{H}^3$, we may assume that $q=0$ and that $D$ is the unit disc 
about the origin. Let $j_{k,n}:D\setminus\{q\}\rightarrow\Bbb{H}^3$ be the solution of the Plateau problem 
$(D\setminus\{q\},z\mapsto z^n)$ with Gaussian curvature equal to $k$. By Lemma $7.2.1$ of \cite{LabA}, 
$(D\setminus\{q\},j_{k,n})$ is a graph over $(\Sigma',i_k)$. In otherwords, there exists a neighbourhood $\Omega$ of $q$ 
in $\Sigma'$, a diffeomorphism $\alpha:\Omega\setminus\{p\}\rightarrow D\setminus\{q\}$ and a smooth function 
$f:\Omega\rightarrow [0,+\infty[$ such that, if $\hat{\mathi}_k$ is the Gauss lifting of $i_k$, then, for all 
$x\in\Omega\setminus\{p\}$:
$$
j_k\circ\alpha(x) = \opExp_x(f(x)\hat{\mathi}_k(x)).
$$
\noindent Using elementary hyperbolic geometry (see, for example, \cite{BallGromSch}), the mapping $\alpha$ is dilating with respect to the metrics induced by the immersions. However,
if $j_{k,1}:D\setminus\{q\}\rightarrow\Bbb{H}^3$ is the solution of the Plateau problem $(D\setminus\{q\},z\mapsto z)$,
then, by uniqueness of solutions, $j_{k,n}$ factors through as an $n$-fold covering of $j_{k,1}$. Finally, by 
Corollary \procref{Control_Of_Area_And_Volume_Of_Rotationally_Symmetric_Cusp}, 
the area of the cusp end of $(D\setminus\{q\},j_{k,1})$ is finite. Thus, the 
area of the cusp of $(\Sigma',i_k)$ about $p$ is finite. The result now follows.\qed
\newhead{Bibliography}
{\leftskip = 5ex \parindent = -5ex
\leavevmode\hbox to 4ex{\hfil\cite{BallGromSch}}\hskip 1ex{Ballman W., Gromov M., Schroeder V., {\sl Manifolds of nonpositive curvature}, Progress in Mathematics, {\bf 61}, Birkh\"auser, Boston, (1985)}
\medskip
\leavevmode\hbox to 4ex{\hfil\cite{Grom}}\hskip 1ex{Gromov M., Foliated plateau problem, part I : Minimal varieties, {\sl GAFA} {\bf 1}, no. 1, (1991), 14--79}%
\medskip
\leavevmode\hbox to 4ex{\hfil\cite{LabC}}\hskip 1ex{Labourie F., Probl\`eme de Minkowski et surfaces \`a 
courbure constante dans les vari\'et\'es hyperboliques, {\sl Bull. Soc. Math. France} {\bf 119}, no. 3, (1991), 307--325}
\medskip
\leavevmode\hbox to 4ex{\hfil\cite{LabB}}\hskip 1ex{Labourie F., Probl\`emes de Monge-Amp\`ere, courbes holomorphes et laminations, {\sl GAFA} {\bf 7}, no. 3, (1997), 496--534}
\medskip
\leavevmode\hbox to 4ex{\hfil \cite{LabA}}\hskip 1ex{Labourie F., Un lemme de Morse pour les surfaces convexes, {\sl Invent. Math.} {\bf 141} (2000), 239--297}
\medskip
\leavevmode\hbox to 4ex{\hfil\cite{RosSpruck}}\hskip 1ex{Rosenberg H., Spruck J. On the existence of convex hyperspheres of constant Gauss curvature in hyperbolic space, {\sl J. Diff. Geom.} {\bf 40} (1994), no. 2,
379--409}
\medskip
\leavevmode\hbox to 4ex{\hfil \cite{SmiB}}\hskip 1ex{Smith G., Hyperbolic Plateau problems, Preprint Univ. Paris XI, 2005,\hfill\break math.DG/0506231}%
\medskip
\leavevmode\hbox to 4ex{\hfil \cite{SmiC}}\hskip 1ex{Smith G., Pointed k-surfaces, to appear in {\it Bull. Soc. Math. France\/}}%
\medskip
\leavevmode\hbox to 4ex{\hfil\cite{SmiE}}\hskip 1ex{Smith G., Th\`ese de doctorat, Paris (2004)}%
\medskip
}%
\enddocument